\documentclass[12pt, reqno]{amsart}
\usepackage{amsmath}
\usepackage{amsthm}
\usepackage{amssymb}
\usepackage{enumerate}
\usepackage{mathrsfs}
\usepackage{graphicx}
\usepackage{graphicx}
\usepackage{amssymb,amsfonts}
\usepackage{amsmath}
\usepackage[colorlinks,linkcolor=blue,anchorcolor=blue,citecolor=blue]{hyperref}
\usepackage[numbers,sort&compress]{natbib}
\usepackage{marvosym}
\usepackage{epstopdf}
\usepackage{caption}
\usepackage{multicol}
\usepackage{multirow}
\usepackage{makebox}
\usepackage{subfigure}
\usepackage{dsfont}
\numberwithin{equation}{section}
\newtheorem{theo}{Theorem} 

\newtheorem{lem}{Lemma}

\begin{document}

\title {$\textrm{SL}(n)$ covariant vector-valued valuations on Orlicz spaces}

\author[Zeng and Lan]{Chunna Zeng$^{*}$ and Yu Lan}

\address{1.School of Mathematical Sciences, Chongqing Normal University, Chongqing 401331, People's Republic of China; Institut f\"{u}r Diskrete Mathematik und Geometrie, Technische Universit\"{a}t Wien, Wiedner Hauptstra\ss e 8--10/1046, 1040 Wien, Austria }
\email{zengchn@163.com}

\address{2.School of Mathematics and Sciences, Chongqing Normal University, Chongqing 401332, People's Republic of China}
\email{18181567221@163.com}

\thanks{The first author is supported by the Major Special Project of NSFC (Grant No. 12141101), the Young Top-Talent program of Chongqing (Grant No. CQYC2021059145), the Science and Technology Research Program of Chongqing Municipal Education Commission (Grant No. KJZD-K202200509) and Natural Science Foundation Project of Chongqing (Grant No. CSTB2024NSCQ-MSX0937)}
\thanks{{\it Keywords}: valuation, $\textrm{SL}(n)$ covariance, Orlicz spaces.}
\thanks{{*}Corresponding author: Chunna Zeng}
\begin{abstract}
A representation theorem for continuous, $\mathrm{SL}(n)$ covariant vector-valued valuations on Orlicz spaces is established. Such valuations are uniquely characterized as moment vectors.
\end{abstract}

\maketitle

\section{Introduction}
 Let $\mathcal{W}$ be a lattice.  A function $\mathrm{Y}: \mathcal{W}\rightarrow \mathbb{R}^n$ is termed a valuation if
 for any  $Q_{1}, Q_{2}$ such that $ Q_{1}\cup Q_{2},$ and $ Q_{1}\cap Q_{2}$ are in $\mathcal{W},$ the following holds:
\begin{equation}
\mathrm{Y}(Q_{1})+\mathrm{Y}(Q_{2})=\mathrm{Y}(Q_{1}\cup Q_{2})+\mathrm{Y}(Q_{1}\cap Q_{2}).
\end{equation}
Valuations on convex bodies originated from Dehn’s solution of  Hilbert’s third problem. A significant development in this field is the Hadwiger characterization theorem, which provides a comprehensive classification of continuous, rigid motion invariant valuations on convex bodies.
This theorem plays a crucial role in both geometric probability theory and integral geometry theory.
It is the basis of studying the results of modern valuation theory. Monika Ludwig has made substantial contributions to the study of body-valued valuations \cite{Ludwig1,Ludwig2,Ludwig3,Ludwig4,Ludwig9}. For a historical development and recent contributions on convex bodies, star-shaped sets, manifolds and so on, see \cite{Alesker1,Alesker2,Alesker3,Alesker4,Klain1,Klain2,Klain3,Klain-Rota,L-Y-Z1,Ludwig5,Ludwig6,Ludwig7,Ludwig8,McMullen1,McMullen2,Ma,Parapatits1,Parapatits2,Wannerer,Wang1,Zeng-Ma,Zeng-Zhou}.
\par  Recently, there has been growing interest in exploring valuations on function spaces.  Assume that $\mathcal{S}$ is a lattice of some functions. For $h_1, h_2\in\mathcal{S},$
the operations $h_1\vee h_2$ and $h_1\wedge h_2$ denote $\max\{h_1, h_2\}$ and $\min\{h_1, h_2\}$, respectively.
A function $\Psi: \mathcal{S}\rightarrow\mathbb{R}^n$ is a valuation if it satisfies
\begin{equation*}\label{0}
\Psi(h_1\vee h_2)+\Psi(h_1\wedge h_2)=\Psi(h_1)+\Psi(h_2)
\end{equation*}
for all $h_1, h_2, h_1\vee h_2, h_1\wedge h_2\in \mathcal{S}.$ Valuations have been extensively studied in various function spaces, including $L^{p}$-spaces \cite{Ludwig6,Tsang1,Tsang2}, Sobolev spaces \cite{Ludwig5,Ludwig7}, the space of Lipschitz functions \cite{Colesanti1,Colesanti2}, spaces of functions of bounded variation \cite{Wang} and Orlicz spaces \cite{Kone}, among others. 
Tsang \cite{Tsang1} was the first to provide a complete classification of continuous, translation invariant real-valued valuation on $L^{p}(\mathbb{R}^n)$, as well as continuous, rotation invariant real-valued valuation on $L^{p}(\mathbb{S}^n)$.  Wang \cite{Wang3}  provided a  classification of vector-valued, continuous and  $\textrm {SL}(n)$ covariant  valuation on $L^{p}(\mathbb{R}^n, \vert x\vert dx).$   Kone  \cite{Kone} extended Tsang's work from $L^p$-spaces to  Orlicz spaces.

\par
 Let $\phi$ be a Young function. The Orlicz space $L^{\phi}(\mathbb{R}^n, |x|dx)$ is defined as the set
  of measurable functions $h:\mathbb{R}^n\rightarrow\mathbb{R}$ such that $\int_{\mathbb{R}^n}\phi(\delta| h(x) |)|x|dx<\infty$ for every $\delta>0$.
An operator $\Psi: L^{\phi}(\mathbb{R}^n,\vert x\vert dx)\rightarrow\mathbb{R}^n$ is termed $\textrm{SL}(n)$ covariant if for all $h\in L^{\phi}(\mathbb{R}^n,\vert x\vert dx)$ and $\vartheta\in\textrm{SL}(n)$,
\begin{equation*}
\Psi(h\circ\vartheta^{-1})=\vartheta \Psi(h).
\end{equation*}
 For a measurable function $h:\mathbb{R}^n\rightarrow\mathbb{R}$, we define the moment vector $m(h)$ as follows:
\begin{equation*}
m(h)=\int_{\mathbb{R}^n} h(x)xdx.
\end{equation*}

In order to investigate function-valued valuations on  Orlicz spaces, it often relies on the results established for convex bodies.
Recently, Zeng and Ma \cite{Zeng-Ma} established the following
$\mathrm{SL}(n)$ covariant vector-valued valuations on polytopes.
\begin{theo}\label{theo1}\cite{Zeng-Ma}
 An operator $\mathrm{Y}: \mathcal{P}^2\rightarrow\mathbb{R}^2$ is a measurable and  $\mathrm{SL(2)}$ covariant valuation if and only if there are real constants $c_{1},c_{2},c_{3},\widetilde{c_{1}},\widetilde{c_{2}},\widetilde{c_{3}}$ such that
\begin{align*}
\mathrm{Y}(Q)
=&c_{1}m(Q)+\widetilde{c_{1}}m([0,Q])+c_{2}e([0,Q])+c_{3}h([0,Q])
+\widetilde{c_{2}}e([0,u_{1},\ldots,u_{r}])\\
&+\widetilde{c_{3}}h([0,u_{1},\ldots,u_{r}])
\end{align*}
for all $Q\in\mathcal{P}^2$, with vertices $u_{1},\ldots,u_{r}$ visible from the origin in counterclockwise order.
\end{theo}
The operator \( e:\mathcal{P}_{0}^{2}\rightarrow\mathbb{R}^2 \) is defined as
follows: \[e(Q) =
\begin{cases}
u + v, & \text{if } \dim Q = 2 \text{ and } Q \text{ has edges } [0, u] \text{ and } [0, v], \ \text{or} \  0\in \text{int} \ P \ \text{with}\\
&  P \ \text{has an edge} \ [u,v]; \\
2(u + v), & \text{if } \dim Q = 1 \text{ and } 0\in Q = [u, v];
\\
0, & \text{otherwise}.
 \end{cases} \]

The operator \( h: \mathcal{P}_{0}^{2} \rightarrow \mathbb{R}^2 \) is defined as follows:

\[
h(Q) =
\begin{cases}
u_1 - u_r, & \text{if } \dim Q = 2, \ Q = [0, u_1, \ldots, u_r] \  \text{with} \\
& \{u_{1}, \ldots, u_{r}\} \text{ arranged counterclockwise} \  \text{and} \ 0 \in \text{bd } Q; \\
0, & \text{if } Q = \{0\}, \text{ a line segment, or } 0 \in \text{int } Q.
\end{cases}
\]

\begin{theo}\label{theo2}\cite{Zeng-Ma}
 An operator $\mathrm{Y}: \mathcal{P}^n\rightarrow\mathbb{R}^n$, where $n\geq3$, is a measurable and $\mathrm{SL}(n)$ covariant valuation if and only if there are real constants $c_{1},c_{2}$ such that
\begin{equation*}
\mathrm{Y}(Q)=c_{1}m(Q)+c_{2}m([0,Q])
\end{equation*}
for every $Q\in\mathcal{P}^n$.
\end{theo}

The purpose of our paper is to characterize continuous, $\textrm{SL}(n)$ covariant, vector-valued valuations on $L^{\phi}(\mathbb{R}^n, \vert x\vert dx)$. Let $C_{0}(\mathbb{R})$ denote the set of continuous functions $\xi:\mathbb{R}\rightarrow\mathbb{R}$ where $\xi(0)=0.$
Similarly, let $C_{\phi}(\mathbb{R})$ represent the set of continuous functions
$\xi:\mathbb{R}\rightarrow\mathbb{R}$ that satisfy $|\xi(\beta)|\leq\lambda\phi(|\beta|)$ for some $\lambda\geq0$ and all $\beta\in\mathbb{R},$  with $\phi$ being a Young function. The detailed definition of the $\Delta_2$ condition will be introduced in Section 2. Unless otherwise specified, we assume that $n\geq 2.$
The following theorem is a vector-valued analog of Kone's characterization of real-valued valuations (see \cite{Kone}).

\begin{theo}\label{theo3}
Suppose $\phi$ is a Young function that satisfies the $\Delta_{2}$ condition.  An operator
$\Psi:L^{\phi}(\mathbb{R}^n,|x|dx)\rightarrow\mathbb{R}^n$
is a continuous, $\mathrm{SL}(n)$ covariant valuation if and only if there is $\xi\in C_{\phi}(\mathbb{R})$ such that
\begin{equation*}
\Psi(h)=m(\xi\circ h)
\end{equation*}
for every $h\in  L^{\phi}(\mathbb{R}^n,|x|dx)$.
\end{theo}


\section{~Preliminaries}
\par 
 In order to investigate the properties of  Orlicz spaces, it is essential to introduce the concept of Young functions.
A Young function is a real-valued, non-negative function $\phi: [0,\infty)\rightarrow[0,\infty)$ that can be expressed as:
\begin{equation*}
\phi(t)=\int_{0}^{t}\varphi(s)ds, \, \,t\geq0,
\end{equation*}
where  $\varphi: [0,\infty)\rightarrow[0,\infty)$ is non-decreasing and satisfies the following conditions:\\
(1) $\varphi(s)>0$ for $s>0$ with $\varphi(0)=0$;\\
(2) $\varphi$ is right continuous for $s\geq0$;\\
(3) $\lim\limits_{s\rightarrow\infty}\varphi(s)=\infty$.\\
The Young function $\phi$ is convex, non-negative, continuous and strictly increasing on $[0,\infty),$ and it satisfies both $\phi(0)=0$ and $\lim\limits_{t\rightarrow\infty}\phi(t)=\infty$. 
The function $\phi^{*}$ is the complementary function to $\phi$, defined by
\begin{equation*}
\phi^{*}(t)=\int_{0}^{t}\psi(s)ds,
\end{equation*}
where $\psi(t)=\mathrm{sup}\{s:\varphi(s)\leq t\}$
for $t\geq0$. It is clear that $\phi^{*}$ is a Young function.
The pair $(\phi,\phi^{*})$
is referred to as  a pair of  complementary Young functions.


\par \par
Given a measure space $(\mathbb{R}^{n},\mathcal{A},\mu)$ and a Young function
$\phi,$ the Orlicz class $ L_{\dagger}^{\phi}(\mathbb{R}^n, |x| dx)$ is defined as
\begin{equation*}
L_{\dagger}^{\phi}(\mathbb{R}^n,|x| dx)=\{h:\mathbb{R}^n\rightarrow\mathbb{R}, \,\mathrm{mesurable}, \, \rho_{\phi}(h)<\infty\},
\end{equation*}
where \begin{equation*}
\rho_{\phi}(h)=\int_{\mathbb{R}^n}\phi(|h(x)|)|x| dx.
\end{equation*}
And the Orlicz space $L^{\phi}(\mathbb{R}^n,|x| dx)$
is defined as
\begin{equation*}
L^{\phi}(\mathbb{R}^n,|x| dx)=\{h:\mathbb{R}^n\rightarrow\mathbb{R},\,\mathrm{mesurable}, \,\rho_{\phi}(\delta h)<\infty\, \,\mathrm{for\,\, every} \,\, \delta>0\}.
\end{equation*}

A Young function $\phi$ is called to satisfy $\Delta_{2}$ condition
if there exist $T\geq0$ and a positive constant $C$ such that for all $t\geq T,$
the function satisfies the inequality:
\begin{equation*}
\phi(2t)\leq C\phi(t).
\end{equation*}
If $\phi$ satisfies $\Delta_{2}$ condition, then the Orlicz class is equivalent to the Orlicz space.

\par 
For $h\in L^{\phi}(\mathbb{R}^n,|x|dx)$, the Orlicz norm $\| h\|_{\phi}$ of $h$ in $L^{\phi}(\mathbb{R}^n,|x| dx)$  is given by
\begin{equation*}
\| h \|_{\phi}=\mathrm{sup}\left\{\int_{\mathbb{R}^n}|h(x)g(x)||x| dx: g\in L^{\phi^{*}}_{\dagger}(\mathbb{R}^n,|x| dx), \rho_{\phi^*}(g)\leq1)\right\}.
\end{equation*}
The functional $\|\cdot\|_{\phi}: L^{\phi}(\mathbb{R}^n,|x|dx)\rightarrow\mathbb{R}$ is a semi-norm. Since we identify $h$ and $g$ when $\| h-g\|_{\phi}=0$, which converts $\| \cdot\|_{\phi}$ into a norm. Thus, $L^{\phi}(\mathbb{R}^n,|x|dx)$ is a normed linear space.

\par 
A sequence of functions $\{h_k\}\subset L^{\phi}(\mathbb{R}^n, |x| dx)$
 is said to converge to a function $h\in L^{\phi}(\mathbb{R}^n, |x| dx)$
 if $\lim\limits_{k\rightarrow \infty} \|h_k - h\|_{\phi} = 0.$
This convergence provides a foundation for analyzing continuity within our function spaces.
Specifically, a function \(\Psi: L^{\phi}(\mathbb{R}^n, |x| dx) \to \mathbb{R}^n\) is considered continuous if, for every sequence $\{h_k\}$ converging to \(h\) in \(L^{\phi}(\mathbb{R}^n, |x| dx)\) as $k\rightarrow\infty$, we have \(\Psi(h_k) \to \Psi(h)\).
\par Next, we introduce the subspace $ L^{\phi}_{a}(\mathbb{R}^n, |x|dx)$ of $L^{\phi}(\mathbb{R}^n, |x|dx) $.
The annulus $A[r,R)$ represents the region between two concentric spheres and is expressed as
\begin{equation*}
A[r,R)=\{x\in\mathbb{R}^n: r\leq |x|<R\},
\end{equation*}
where $r$ and $R$ are real numbers satisfying $0\leq r<R,$ and $|x|$
 is the distance from the point $x$ to the origin.
 The subspace \( L^{\phi}_{a}(\mathbb{R}^n, |x|dx) \) is defined by
\[
L^{\phi}_{a}(\mathbb{R}^n, |x|dx) = \{ h \in L^{\phi}(\mathbb{R}^n, |x|dx) : \text{supp}(h) \subseteq A[r, R) \},
\]
where 
\[
\text{supp}(h) = \{ x \in \mathbb{R}^n : h(x) \neq 0 \}.
\]
In this context, the condition \( \chi_{\mathbb{R}^n \setminus A[r, R)} h = 0 \) implies that \( h \) vanishes outside the annulus \( A[r, R) \), thereby ensuring that \( h \) is supported within this region.
\par In this paper, we give two measurable spaces: $(\mathbb{R}^n,\mathscr{B},\lambda)$ and $(\mathbb{R}^n,\mathscr{B},\mu_{n})$. Here, $\mathscr{B}$ represents the Lebesgue measurable sets, and $\lambda$ denotes the Lebesgue measure. The measure $\mu_{n}$ is defined for any set $M\in\mathscr{B}$ as
\begin{equation*}
\mu_{n}(M)=\int_{M}\vert x \vert dx.
\end{equation*}
We write $L^{\phi}(\mathbb{R}^n)$ as $L^{\phi}(\lambda)$ and use $L^{\phi}(\mu_{n})$ to denote $L^{\phi}(\mathbb{R}^n, |x| dx)$. Thus, for
  $h\in L^{\phi}(\mathbb{R}^n)$, we typically write $\int_{\mathbb{R}^n}h(x)dx$ instead of $\int_{\mathbb{R}^n}hd\lambda$. Similarly, 
  we replace $\int_{\mathbb{R}^n}h d\mu_{n} $ by $\int_{\mathbb{R}^n}h(x)|x|dx$. For convenience, we refer to both $\lambda$-measurable functions and $\mu_{n}$-measurable functions as measurable functions. It is noteworthy that if $M\in\mathscr{B}$, then $\mu_{n}(M)=0$  is equivalent to $\lambda(M)=0$. Thus, we abbreviate almost everywhere with respect to $\mu_{n}$ and  $\lambda$ as a.e.

\par
The indicator functions and simple functions play crucial roles in our subsequent proofs. The indicator function \(\chi_M\) for a set \(M \subset \mathbb{R}^n\) is defined as follows
\[
\chi_M(x) =
\begin{cases}
1, & \text{if } x \in M, \\
0, & \text{otherwise}.
\end{cases}
\]
If $s$ is a simple function,  there exist pairwise disjoint sets $M_{1},\ldots,M_{l}\in\mathscr{B}$ that have finite measure and scalars $\alpha_{i}$ for $i=1,\ldots,l$ such that
\begin{equation*}
s=\sum\limits_{i=1}\limits^{l}\alpha_{i}\chi_{M_{i}}(x).
\end{equation*}

\par Denote
$e_{1},\ldots,e_{n}$ as the standard basis vectors in $\mathbb{R}^n$. Let $[u_{1},\ldots,u_{k}]$ represent the convex hull of $u_{1},\ldots,u_{k}\in\mathbb{R}^n$. Define $\mathcal{P}^n$ as the set of compact convex polytopes, and $\mathcal{P}_{0}^{n}$ as its subset consisting of compact convex polytopes that contain the origin 0, both equipped with the topology derived from the Hausdorff metric. Denote $m(Q)$ as moment vector if
 \begin{equation*}
m(Q)=\int_{Q}xdx
 \end{equation*}
 for $Q\in\mathcal{P}^{n}.$

 An operator $\mathrm{Y}:\mathcal{P}^n\rightarrow \mathbb{R}^{n}$ is $\mathrm{SL}(n)$ covariant if it satisfies
\begin{equation*}
\mathrm{Y}(\vartheta Q)=\vartheta \mathrm{Y}(Q)
\end{equation*}
for any $Q\in \mathcal{P}^n$ and $\vartheta\in \mathrm{SL}(n)$. An operator $\mathrm{Y}:\mathcal{P}^n\rightarrow\mathbb{R}^n$ is termed simple if $\mathrm{Y}(Q)=0$ for any $Q\in\mathcal{P}^n$, where  dim $Q$ is at most $n-1$. It is well known that the moment vector $m$ is simple, $\mathrm{SL}(n)$ covariant and continuous function.

\par
The following quotes play an important role in our proofs.

\begin{lem}\label{lem1}\cite{Kone}
If $\phi$ is a Young function, then
\begin{equation*}
\lim\limits_{t\rightarrow\infty}\frac{\phi(t)}{t}=\infty \,\,and \,\,\lim\limits_{t\rightarrow\infty}\frac{\phi^{-1}(t)}{t}=0
\end{equation*}
for all $t\geq0$.
\end{lem}

\par A closed cube $C$ in $\mathbb{R}^n$ is defined as
\begin{equation*}
C=\{(x_{1},\ldots,x_{n}): c_{i}\leq x_{i}\leq d_{i},\  x_{i}, c_{i}, d_{i}\in\mathbb{R}\},
\end{equation*}
where $\vert d_{i}-c_{i}\vert= \vert c_{j}-d_{j}\vert$,  $1\leq i, j\leq n$. For subsets $Q_1,Q_2\subseteq\mathbb{R}^n$, the definition of the symmetric difference of  $Q_1$ and $ Q_2$ is
 \begin{equation*}
 Q_1\bigtriangleup Q_2=Q_1\setminus Q_2 \cup Q_2\setminus Q_1.
 \end{equation*}

\begin{lem}\label{lem2}\cite{Tsang1}
Suppose that $\varepsilon>0$ and $D\in\mathscr{B}$ with $\lambda(D)<\infty$. Then there exists a series of closed cubes $C_{1}, C_{2},\ldots, C_{l}\subseteq\mathbb{R}^n$ such that $\lambda(D\bigtriangleup \bigcup\limits_{i=1}\limits^{l} C_{i})<\varepsilon$.
\end{lem}

\begin{lem}\label{lem3}\cite{Kone}
Let $M\subset\mathbb{R}^n$. If $0<\mu_{n}(M)<\infty$, then
\begin{equation*}
\|\chi_{M}\|_{\phi}=\mu_{n}(M)\phi^{*^{-1}}\left(\frac{1}{\mu_{n}(M)}\right).
\end{equation*}
 \end{lem}


\begin{lem}\label{lem6}\cite{Kone}
 Let $\{h_{i}\}\subset L^{\phi}(\mu_{n})$ be a sequence and $h\in L^{\phi}(\mu_{n})$ such that $\Vert h_{i}-h\Vert_{\phi}\rightarrow0$. Then, there are a subsequence $\{h_{i_{j}}\}$ of $\{h_{i}\}$ and a function $g^*\in L^{\phi}(\mu_{n})$ satisfying the following conditions:

 $(1)\ h_{i_{j}}\rightarrow h$ a.e. for all $j$;

$(2)\ |h_{i_{j}}|\leq g^*$ a.e.
\end{lem}

\section{Proof of main theorems}
\begin{theo}\label{lem7} Let $\xi \in C_{0}(\mathbb{R})$. Assume that for every $\delta>0$, there is a real number $\lambda_{\delta}\geq 0$  such that $|\xi(\beta)|\leq\lambda_{\delta}\phi(\delta|\beta|)$ for any $\beta \in \mathbb{R}$. Then the operator $\Psi:L^{\phi}(\mu_{n})\rightarrow \mathbb{R}^n$ defined by
\begin{equation*}
\Psi(h)=m(\xi\circ h)
\end{equation*}
is a continuous and $\mathrm{SL}(n)$ covariant valuation on $L^{\phi}(\mu_{n})$.
\end{theo}
\emph{Proof.}
 Since there is a real number $\lambda_{\delta}\geq0$ such that $| \xi(\beta)|\leq\lambda_{\delta}\phi(\delta|\beta|)$ for all $\beta\in\mathbb{R}$, we have
\begin{align*}
|m(\xi\circ h)|&=\left|\int_{\mathbb{R}^{n}}(\xi\circ h)(x)xdx\right|\leq \int_{\mathbb{R}^{n}}|(\xi\circ h)(x)||x|dx\\
&\leq \lambda_{\delta}\int_{\mathbb{R}^{n}}\phi(\delta|h(x)|)|x|dx<\infty.
\end{align*}
Therefore, it follows that $|\Psi(h)|<\infty$.  From $\xi(0)=0$, we have $\Psi(0)=\mathrm{\textbf{0}}.$
 Let $h_{1},h_{2}\in L^{\phi}(\mu_{n})$, then
\begin{align*}
\Psi(h_{1}\vee h_{2})+\Psi(h_{1}\wedge h_{2})&=\int_{\mathbb{R}^n}\xi\circ(h_{1}\vee h_{2})x dx+\int_{\mathbb{R}^n}\xi\circ(h_{1}\wedge h_{2})xdx\\
&=\int_{\{h_{1}\geq h_{2}\}}(\xi\circ h_{1})(x)x dx+\int_{\{h_{1}<h_{2}\}}(\xi\circ h_{2})(x) xdx\\
&\ \ \ \ +\int_{\{h_{1}\geq h_{2}\}}(\xi\circ h_{2})(x)xdx+\int_{\{h_{1}<h_{2}\}}(\xi\circ h_{1})(x)x dx\\
&=\int_{\mathbb{R}^n}(\xi\circ h_{1})x dx+\int_{\mathbb{R}^n}(\xi\circ h_{2})(x) xdx\\
&=\Psi(h_{1})+\Psi(h_{2}).
\end{align*}
Thus, $\Psi$ is a valuation on $L^{\phi}(\mu_{n})$.

 We now aim to show that $\Psi$ is $\mathrm{SL}(n)$ covariant. For any $h\in L^{\phi}(\mu_{n})$ and  $\vartheta\in \mathrm{SL}(n),$ we have
\begin{equation*}
\Psi(h\circ\vartheta^{-1})=\int_{\mathbb{R}^n}(\xi\circ h)(\vartheta^{-1}x)xdx=\vartheta\int_{\mathbb{R}^n}(\xi\circ h)(x)xdx=\vartheta \Psi(h).
\end{equation*}

\par
Suppose that $\{h_{i}\}$ is a sequence of operators in $L^{\phi}(\mu_{n}),$ and $h_{i}\rightarrow h$ as $i\rightarrow\infty$ in $L^{\phi}(\mu_{n})$. To establish the continuity of $\Psi$, we need to demonstrate that any subsequence $\{\Psi(h_{i_{j}})\}$ of $\{\Psi(h_{i})\} $ has a subsequence $\{\Psi(h_{i_{j_{k}}})\}$ converges to $\{\Psi(h)\}$. Since \{$h_{i_{j}}\}$ is a subsequence of $\{h_{i}\}$, we also have $h_{i_{j}}\rightarrow h$ in $L^{\phi}(\mu_{n})$. By Lemma \ref{lem6}, there is a subsequence $\{h_{i_{j_{k}}}\}$ of $\{h_{i_{j}}\}$ satisfying that $h_{i_{j_{k}}}\rightarrow h$ a.e., and $g^*\in L^{\phi}(\mu_{n})$ such that $|h_{i_{j_{k}}}|\leq g^*$ a.e. for all $k$. Since $\xi$ is continuous, it follows that
\begin{equation*}
(\xi\circ h_{i_{j_{k}}})(x)x\rightarrow(\xi\circ h) (x)x \ \  a.e.
\end{equation*}
Additionally,
because  $|\xi(\beta)|\leq\lambda_{\delta}\phi(\delta\vert \beta\vert)$ for all $\beta\in \mathbb R$, we obtain
\begin{equation*}
|(\xi\circ h_{i_{j_{k}}})(x)x|\leq|(\xi\circ h_{i_{j_{k}}})(x)||x|
\leq\lambda_{\delta}\phi(\delta|h_{i_{j_{k}}}(x)|)|x|\leq\lambda_{\delta}\phi(\delta g^*)|x|,\, a.e.
\end{equation*}
Since $\lambda_{\delta}\phi(\delta g^*)|x|$ is integrable, by the dominated convergence theorem, we obtain
\begin{equation*}
\lim\limits_{k\rightarrow\infty}\int_{\mathbb{R}^n}(\xi\circ h_{i_{j_{k}}})(x)xdx=\int_{\mathbb{R}^n}(\xi\circ h)(x)xdx.
\end{equation*}
Thus, $\Psi(h_{i_{j_{k}}})\rightarrow \Psi(h)$ as $h_{i_{j_{k}}}\rightarrow h$ in $L^{\phi}(\mu_{n})$, which implies that  $\Psi$ is continuous.
\qed
\begin{lem}\label{lem8}
An operator $\mathrm{Y}:\mathcal{P}^n\rightarrow\mathbb{R}^n$ is a  simple, continuous and $\mathrm{SL}(n)$ covariant valuation if and only if there is a constant $c\in\mathbb{R}$ such that
\begin{equation*}
\mathrm{Y}(Q)=cm(Q)
\end{equation*}
for all $Q\in \mathcal{P}^n$.
\end{lem}
\emph{Proof.}
Let $Q,Q_{i}\in\mathcal{P}^n$ with $Q_{i}\rightarrow Q$ as $i\rightarrow\infty$. Observe that $[0,Q_{i}]\rightarrow [0,Q]$ as $i\rightarrow\infty$. 
 Theorem \ref{theo1} and Theorem \ref{theo2} establish a classification of measurable and $\mathrm{SL}(n)$ covariant valuation on $\mathcal{P}^n$. In the case of $n=2$, if the operator  $\mathrm{Y}:\mathcal{P}^2\rightarrow\mathbb{R}^2$ is a measurable and $\mathrm{SL}(2)$ covariant valuation, then there are real constants $c_{1},  c_{2}, c_{3},\widetilde{c_{1}}, \widetilde{c_{2}}, \widetilde{c_{3}}$ such that
\begin{align*}
\mathrm{Y}(Q)&=
c_{1}m(Q)+\widetilde{c_{1}}m([0,Q])+c_{2}e([0,Q])+\widetilde{c_{2}}e([0,u_{1},\ldots,u_{r}])
+c_{3}h([0,Q])\\
&+\widetilde{c_{3}}h([0,u_{1},\ldots,u_{r}])
\end{align*}
for every $Q\in\mathcal{P}^2$, where the vertices $u_{1},\ldots,u_{r}$ are visible from the origin and labeled counterclockwise.
 Let
\begin{equation}\label{09}
\widetilde{Y}(Q)=c_{2}e([0,Q])+\widetilde{c_{2}}e([0,u_{1},\ldots,u_{r}])+c_{3}h([0,Q])+\widetilde{c_{3}}h([0,u_{1},\ldots,u_{r}])
\end{equation}
for $Q\in\mathcal{P}^2$. Since $\mathrm{Y}$ is continuous on $\mathcal{P}^2$ and $m$ is also continuous, it leads that $\widetilde{Y}$ is a continuous function.
\par Consider $Q=[0, e_{1}]$ and $Q_{\varepsilon}=[e_{1},\varepsilon e_{2}]$, where $Q_{\varepsilon}\rightarrow Q$ as $\varepsilon\rightarrow0$. Substituting $Q$ and $Q_{\varepsilon}$ into (\ref{09}), we have
\begin{equation*}
\widetilde{Y}(Q)=2c_{2}e_{1},
\end{equation*}
and
\begin{equation*}
\widetilde{Y}(Q_{\varepsilon})=(c_{2}+\widetilde{c_{2}}+c_{3}+\widetilde{c_{3}})e_{1}+(c_{2}+\widetilde{c_{2}}-c_{3}-\widetilde{c_{3}})\varepsilon e_{2}.
\end{equation*}
Due to the continuity of $\widetilde{Y}$, it follows that $\widetilde{Y}(Q_{\varepsilon})\rightarrow\widetilde{Y}(Q)$ as $\varepsilon\rightarrow0$, which gives
\begin{equation*}
c_{2}+\widetilde{c_{2}}+c_{3}+\widetilde{c_{3}}=2c_{2},
\end{equation*}
leading to
\begin{equation}\label{19}
\widetilde{c_{2}}+c_{3}+\widetilde{c_{3}}=c_{2}.
\end{equation}
 Similarly, we take $Q'=[0, e_{2}]$ and $Q'_{\varepsilon}=[\varepsilon e_{1}, e_{2}]$, where $Q'_{\varepsilon}\rightarrow Q'$ as $\varepsilon\rightarrow0$,  substituting $Q'$ and $Q'_{\varepsilon}$  into (\ref{09}). Then
 \begin{equation*}
\widetilde{Y}(Q')=2c_{2}e_{2},
\end{equation*}
and
\begin{equation*}
\widetilde{Y}(Q'_{\varepsilon})=(c_{2}+\widetilde{c_{2}}+c_{3}+\widetilde{c_{3}})\varepsilon e_{1}+(c_{2}+\widetilde{c_{2}}-c_{3}-\widetilde{c_{3}})e_{2}.
\end{equation*}
From the continuity of $\widetilde{Y}$, so  $\widetilde{Y}(Q'_{\varepsilon})\rightarrow\widetilde{Y}(Q')$ as $\varepsilon \rightarrow 0$. We obtain
\begin{equation*}
c_{2}+\widetilde{c_{2}}-c_{3}-\widetilde{c_{3}}=2c_{2},
\end{equation*}
which leads to
\begin{equation}\label{29}
\widetilde{c_{2}}-c_{3}-\widetilde{c_{3}}=c_{2}.
\end{equation}
Next, consider $Q''_{\varepsilon}=[e_{1},e_{2},-\varepsilon e_{1}]$ and $Q'''_{\varepsilon}=[\varepsilon e_{1}, e_{2},- e_{1}]$, where $Q''_{\varepsilon}\rightarrow Q''=[0,e_{1},e_{2}] $ and $Q'''_{\varepsilon}\rightarrow Q'''=[0,e_{2},-e_{1}]$  as $\varepsilon\rightarrow0$, respectively. Similarly, we have
 \begin{equation}\label{39}
c_{2}-c_{3}=0,
\end{equation}
and
\begin{equation}\label{49}
c_{2}+c_{3}=0.
\end{equation}
From (\ref{19}), (\ref{29}), (\ref{39}) and (\ref{49}), it follows that
\begin{equation*}
c_{2}=\widetilde{c_{2}}=c_{3}=\widetilde{c_{3}}=0.
\end{equation*}
Thus,  $\mathrm{Y}(Q)=c_{1}m(Q)+\widetilde{c_{1}}m([0,Q])$.
 Let $Q=[e_{1},e_{2}]\in\mathcal{P}^2$,
 then
\begin{equation*}
\mathrm{Y}(Q)=\widetilde{c_{1}}m([0,Q])=\widetilde{c_{1}}\left(\frac{1}{6},   \frac{1}{6}\right).
\end{equation*}
Notice that if $\mathrm{Y}$ is a simple  valuation, it implies that
 $\widetilde{c_{1}}=0.$

If $n\geq3$, $\mathrm{Y}:\mathcal{P}^n\rightarrow\mathbb{R}^n$ is a measurable and $\mathrm{SL}(n)$ covariant valuation, then there exist constants $c_{1},  c_{2}\in\mathbb{R}$ such that
\begin{equation}\label{69}
\mathrm{Y}(Q)=c_{1}m(Q)+c_{2}m([0,Q]).
\end{equation}
Now, let $Q=[e_{1},\ldots,e_{n}]$ in (\ref{69}), and observe that dim $Q=n-1$. Since $m$ is a simple function, we have
\begin{equation*}
\mathrm{Y}(Q)=c_{2}m([0,Q])=c_{2}\left(\frac{1}{(n+1)!},\frac{1}{(n+1)!}\right).
\end{equation*}
Due to the simplicity of $\mathrm{Y}$, it leads to  $c_{2}=0$.
Therefore, $\mathrm{Y}(Q)=c_{1}m(Q)$.
\qed


\begin{lem}\label{lem4}
 An operator $\Psi:L^{\phi}(\mu_{n})\rightarrow\mathbb{R}^n$ is an $\mathrm{SL}(n)$ covariant valuation, then
\begin{equation*}
\Psi(0)=\mathrm{\textbf{0}}.
\end{equation*}
\end{lem}
\emph{Proof.} For a real number $k\neq0$, let
\begin{equation*}
\vartheta_{n}=
\begin{pmatrix}
k &   & &\\
  & k &  &\\
  &   & \ddots &\\
  &   &  & k^{1-n}
\end{pmatrix},
\end{equation*}
which belongs to  $\mathrm{SL}(n)$. Because $\Psi$ is $\mathrm{SL}(n)$ covariant, it follows that
\begin{equation*}
\Psi(0\circ\vartheta_{n}^{-1})=\Psi(0)=\vartheta_{n}\Psi(0).
\end{equation*}
Let $\Psi(0)=(y_{1},\ldots,y_{n})^{t}$. Then, we have
\begin{equation*}
(y_{1},\ldots,y_{n-1},y_{n})^{t}=(ky_{1},\ldots,ky_{n-1},k^{1-n}y_{n})^{t}.
\end{equation*}
Since $k\neq0$, this implies that $y_{i}=0$ for $1\leq i\leq n$.
\begin{lem}\label{lem10}
If an operator $\Psi:L^{\phi}(\mu_{n})\rightarrow\mathbb{R}^n$ is a continuous and $\mathrm{SL}(n)$ covariant valuation, then there is $\xi\in C_{0}(\mathbb{R})$
 such that for any $\alpha\in\mathbb{R}$,
\begin{equation*}
\Psi(\alpha\chi_{Q})=\xi(\alpha)m(Q)
\end{equation*}
for all $Q\in\mathcal{P}^n$.
\end{lem}
\emph{Proof.}
For  $\alpha\in\mathbb{R}$, define $\mathrm{Y}_{\alpha}:\mathcal{P}^n\rightarrow\mathbb{R}^n$ by setting
\begin{equation*}
\mathrm{Y}_{\alpha}(Q)=\Psi(\alpha\chi_{Q})
\end{equation*}
for $Q\in \mathcal{P}^n.$
Since $\Psi$ is a valuation on $L^{\phi}(\mu_{n})$, for $Q_{1},Q_{2},Q_{1}\cup Q_{2},Q_{1}\cap Q_{2}\in\mathcal{P}^n$, we have
\begin{align*}
\mathrm{Y}_{\alpha}(Q_{1}\cup Q_{2})+\mathrm{Y}_{\alpha}(Q_{1}\cap Q_{2}) &=\Psi(\alpha\chi_{Q_{1}\cup Q_{2}})+\Psi(\alpha\chi_{Q_{1}\cap Q_{2}})\\
&=\Psi(\alpha\chi_{Q_{1}}\vee \alpha\chi_{Q_{2}})+\Psi(\alpha\chi_{Q_{1}}\wedge \alpha\chi_{Q_{2}})\\
&=\mathrm{Y}_{\alpha}(Q_{1})+\mathrm{Y}_{\alpha}(Q_{2}).
\end{align*}
Thus, $\mathrm{Y}_{\alpha}:\mathcal{P}^n\rightarrow \mathbb{R}^n$ is a valuation. For any  $Q\in\mathcal{P}^n$ and $\vartheta\in\mathrm{SL}(n)$, we have
\begin{equation*}
\mathrm{Y}_{\alpha}(\vartheta Q)=\Psi(\alpha\chi_{\vartheta Q})=\Psi((\alpha\chi_{Q})\circ\vartheta^{-1})=\vartheta\Psi(\alpha\chi_{Q})
=\vartheta\mathrm{\mathrm{Y_{\alpha}}}(Q).
\end{equation*}
 Therefore, $\mathrm{Y}_{\alpha}$ is $\mathrm{SL}(n)$ covariant.
For $Q\in\mathcal{P}^n$ with dim $Q\leq n-1$, we have $\alpha\chi_{Q}=0$ a.e. in $L^{\phi}(\mu_{n})$. So
\begin{equation*}
\mathrm{Y}_{\alpha}(Q)=\Psi(\alpha\chi_{Q})=\Psi(0)=\textbf{0}.
\end{equation*}
Hence $\mathrm{Y}_{\alpha}$ is simple.
Next, we prove that $\mathrm{Y}_{\alpha}$ is continuous. According to Lemma \ref{lem2}, for $Q\in\mathcal{P}^n$, there exists a series of  $\{C_{i}\}$, a  union of finite closed cubes, such that $\lambda(C_{i}\bigtriangleup Q)<{1}/{i}$ for $i\in\mathbb{N}$. If $x\in C_{i}\bigtriangleup Q$,  there is a positive constant $a$ such that $|x|\leq a.$  Observe that
\begin{equation*}
\mu_{n}(C_{i}\bigtriangleup Q)=\int_{C_{i}\bigtriangleup Q}|x|dx\leq a\int_{C_{i}\bigtriangleup Q}dx= a\lambda(C_{i}\bigtriangleup Q),
\end{equation*}
thus we have $\mu_{n}(C_{i}\bigtriangleup Q)<a/i$. If $0<\mu_{n}(C_{i}\bigtriangleup Q)<\infty$, we obtain that
\begin{equation*}
\|\alpha\chi_{C_{i}}-\alpha\chi_{Q}\|_{\phi}
=|\alpha| \|\chi_{C_{i}\bigtriangleup Q}\|_{\phi}
=|\alpha| \mu_{n}(C_{i}\bigtriangleup Q)\phi^{*^{-1}}\left(\frac{1}{\mu_{n}(C_{i}\bigtriangleup Q)}\right).
\end{equation*}
 Let $t=1/\mu_{n}(C_{i}\bigtriangleup Q)$. By Lemma \ref{lem1}, it follows that $\|\alpha\chi_{C_{i}}-\alpha\chi_{Q}\|_{\phi}\rightarrow0$
as $i\rightarrow\infty$. On the other hand, if $\mu_{n}(C_{i}\bigtriangleup Q)=0$, then $\|\chi_{C_{i}\bigtriangleup Q}\|_{\phi}=0$. Thus, $\|\alpha\chi_{C_{i}}-\alpha\chi_{Q}\|_{\phi}=0$. In both cases, the continuity of $\Psi$ implies that $\mathrm{Y}_{\alpha}$ is continuous on $\mathcal{P}^n$.

\par By Lemma \ref{lem8}, there is a real constant  $c_{\alpha}$ such that
\begin{equation}\label{59}
\Psi(\alpha\chi_{Q})=c_{\alpha}m(Q)
\end{equation}
for $Q\in\mathcal{P}^n$. We define the operator $\xi:\mathbb{R}\rightarrow\mathbb{R}$ by
\begin{equation*}
\xi(\alpha)=c_{\alpha}.
\end{equation*}
Lemma \ref{lem4} shows that $\Psi(0)=\textbf{0}$. In (\ref{59}), we assume that $\alpha=0$ and $m(Q)\neq0$ for $Q\in\mathcal{P}^n$, then it follows that $\xi(0)=0$.
Now we show that the continuity of $\xi$. Take $\alpha\in\mathbb{R}$ and a sequence $\{\alpha_{i}\}\subset \mathbb{R}$ satisfying $\alpha_{i}\rightarrow\alpha$ as $i\rightarrow\infty$. Let $Q\in\mathcal{P}^n$ with  $\lambda(Q)\neq0$, we also have $\mu_{n}(Q)\neq0$.
Observe that
\begin{equation*}
\|\alpha_{i}\chi_{Q}-\alpha\chi_{Q}\|_{\phi}
=|\alpha_{i}-\alpha|\mu_{n}(Q)\phi^{*^{-1}}\left(\frac{1}{\mu_{n}(Q)}\right).
 \end{equation*}
So $\|\alpha_{i}\chi_{Q}-\alpha\chi_{Q}\|_{\phi}\rightarrow0$
 as $\alpha_{i}\rightarrow\alpha$. By the continuity of $\Psi$, it follows that $\Psi(\alpha_{i}\chi_{Q})\rightarrow\Psi(\alpha\chi_{Q})$, which implies $\xi(\alpha_{i})\rightarrow\xi(\alpha)$ as $i\rightarrow\infty$.
\qed\\

The following lemma closely aligns with Wang’s result (\cite{Wang3}, Lemma 3.3), so the proof is omitted.
\begin{lem}\label{lem12}
If  $\Psi:L^{\phi}(\mu_{n})\rightarrow\mathbb{R}^n$ is a valuation, then for   all $\alpha_{1},\ldots,\alpha_{l}\in\mathbb{R},$ $l\in\mathbb{N},$ and  pairwise disjoint sets $M_{1},\ldots,M_{l}\in\mathscr{B}$ with finite measure,
\begin{equation*}
\Psi\left(\sum\limits_{i=1}\limits^{l}\alpha_{i}\chi_{M_{i}}\right)
=\sum\limits_{i=1}\limits^{l}\Psi(\alpha_{i}\chi_{M_{i}}).
\end{equation*}
\end{lem}

\begin{lem}\label{lem13}
If $\Psi: L^{\phi}(\mu_{n})\rightarrow\mathbb{R}^n$ is a continuous and $\mathrm{SL}(n)$ covariant valuation, then there is  $\xi\in C_{0}(\mathbb{R})$ such that for any $\alpha\in\mathbb{R}$,
\begin{equation*}
\Psi(\alpha\chi_{M})=\xi(\alpha)m(M)
\end{equation*}
for $M\in\mathscr{B}$ with $M\subset A[0,R)$.
\end{lem}

\emph{Proof.}
Let  $M\in\mathscr{B}$ with $M\subset A[0,R)$, we have $\lambda(M)<\infty$. According to Lemma \ref{lem2}, there is a series of  $\{K_{i}\}$, which is a  union of finite closed cubes satisfying $\lambda(K_{i}\bigtriangleup M)<{1}/{i}$. If $x\in K_{i}\bigtriangleup M$, there is a constant $a>0$ satisfying that $\vert x\vert\leq a$. Thus we have $\mu_{n}(K_{i}\bigtriangleup M)\leq a\lambda(K_{i}\bigtriangleup M)<{a}/{i}$. If $0<\mu_{n}(K_{i}\bigtriangleup M)<\infty$, note that
\begin{equation*}
\|\alpha\chi_{K_{i}}-\alpha\chi_{M}\|_{\phi}
=|\alpha|\|\chi_{K_{i}\bigtriangleup M}\|_{\phi}=|\alpha|\mu_{n}(K_{i}\bigtriangleup M)\phi^{*^{-1}}\left(\frac{1}{\mu_{n}(K_{i}\bigtriangleup M)}\right).
\end{equation*}
Let $t=1/\mu_{n}(K_{i}\bigtriangleup M)$. By Lemma \ref{lem1}, we have $\alpha\chi_{K_{i}}\rightarrow\alpha\chi_{M}$ in $L^{\phi}(\mu_{n})$.
 If  $\mu_{n}(K_{i}\bigtriangleup M)=0$, we obtain the same result.
Additionally, for any $K_{i}$, there exist pairwise disjoint interiors sets $Q_{i_{1}},\ldots,Q_{i_{l}}\in\mathcal{P}^n$  such that $K_{i}=\bigcup\limits_{j=1}\limits^{l}Q_{i_{j}}$.  Observe that $\alpha\chi_{Q\setminus\mathring{Q}}=0$ a.e. in $L^{\phi}(\mu_{n})$. Hence,
\begin{equation*}
\alpha\chi_{K_{i}}=\sum\limits_{j=1}^{l}\alpha\chi_{\mathring{Q}_{i_{j}}}=\sum\limits_{j=1}^{l}\alpha\chi_{Q_{i_{j}}} \ \ \textrm{a.e.},
\end{equation*}
where $\mathring{Q}_{i_{j}}$ denote the interior of $Q_{i_{j}}$. 
By Lemma \ref{lem4} , Lemma \ref{lem10} and Lemma \ref{lem12}, we obtain
\begin{align*}
\Psi(\alpha\chi_{K_{i}})&=\sum\limits_{j=1}\limits^{l}\Psi(\alpha\chi_{\mathring{Q}_{i_{j}}})\\
&=\sum\limits_{j=1}\limits^{l}\Psi(\alpha\chi_{Q_{i_{j}}})\\
&=\sum\limits_{j=1}\limits^{l}\xi(\alpha)\int_{Q_{i_{j}}}xdx
=\xi(\alpha)m(K_{i}).
\end{align*}
Meanwhile,
\begin{align*}
|\Psi(\alpha\chi_{K_{i}})-\xi(\alpha)m(M)|
&=|\xi(\alpha)||m(K_{i})-m(M)|\\
&\leq|\xi(\alpha)|\int_{\mathbb{R}^n}|(\chi_{K_{i}}-\chi_{M})(x)||x| dx\\
&\leq a\left|\xi(\alpha)\right|\lambda(K_{i}\triangle M).
\end{align*}
 Therefore, $\Psi(\alpha\chi_{K_{i}})\rightarrow \xi(\alpha)m(M)$ as $i\rightarrow\infty$. Due to the continuity of $\Psi$, it follows that
\begin{equation*}
\lim\limits_{i\rightarrow\infty}\Psi(\alpha\chi_{K_{i}})=\xi(\alpha)m(M)=\Psi(\alpha\chi_{M}).
\end{equation*}
\qed



 \par Denote $B_{n}$ as the unit ball in $\mathbb{R}^n$.
 Tsang (\cite{Tsang2}) constructed a collection of pairwise disjoint balls defined by
 \begin{equation}\label{79}
 \{M_{j}:M_{j}=r_{j}B_{n}+c_{j}e_{1}\}
 \end{equation}
 in $\mathbb{R}^n$, satisfying the following conditions: $(1)\ c_{j}>r_{j}$ for $j\geq1$; $(2)\ 0<r_{j}<1$ for $j\geq1$; $(3)\  c_{j}>c_{j-1}+2 $ for $j\geq2$;\
 $(4)\ (c_{j}+r_{j})\lambda(M_{j})=\beta'_{j}$ for $j\geq1$, where $\{\beta_{j}'\}$ is a sequence of positive numbers in $ \mathbb{R}; (5)\ \lim\limits_{j\rightarrow\infty}\frac{c_{j}}{c_{j}+r_{j}}=1$.

 \begin{lem}\label{lem15}
Suppose that $\phi$ is a Young function satisfying the $\Delta_{2}$ condition. If $\xi\in C_{0}(\mathbb R)$ and $|m(\xi\circ h)|<\infty$ for every $h\in L^{\phi}(\mu_{n})$,
then $\xi\in C_{\phi}(\mathbb{R})$.
\end{lem}
\emph{Proof.}
 Assume that for any $\lambda\geq0$, there exists some $\beta_{0}\neq0$ such that $|\xi(\beta_{0})|>\lambda\phi(|\beta_{0}|)$. Consequently, there is a  sequence $\{\beta_{i}\}$ in $\mathbb{R}\setminus\{0\}$  satisfying
\begin{equation*}
|\xi(\beta_{i})|>2^{i}\phi(|\beta_{i}|).
\end{equation*}
This implies that there exists a subsequence $\{\beta_{i_{j}}\}$ of $\{\beta_{i}\}$ satisfying either
$$\xi(\beta_{i_{j}})>2^{i_{j}}\phi(|\beta_{i_{j}}|)
\quad \text{or} \quad   \xi(\beta_{i_{j}})<-2^{i_{j}}\phi(| \beta_{i_{j}}|).$$
In the first case, when
 $\xi(\beta_{i_{j}})>2^{i_{j}}\phi(|\beta_{i_{j}}|)$ holds, we
 define a sequence of pairwise disjoint balls based on the conditions in  (\ref{79}) and  $\beta'_{j}=\frac{1}{2^{i_{j}}\phi(|\beta_{i_{j}}|)}$ for all $j\in\mathbb{N}$.
We define $h:\mathbb{R}^n\rightarrow\mathbb{R}$ by
\begin{equation*}
h=\sum\limits_{j=1}\limits^{\infty}\beta_{i_{j}}\chi_{M_{j}}.
\end{equation*}
 Next, we compute
\begin{equation*}
\begin{aligned}
\int_{\mathbb{R}^n}\phi(|h(x)|)|x| dx
&=\int_{\mathbb{R}^n}\phi\left(\left|\sum\limits_{j=1}\limits^{\infty}
\beta_{i_{j}}\chi_{M_{j}}\right|\right)|x|dx\\
&=\sum\limits_{j=1}\limits^{\infty}\phi(|\beta_{i_{j}}|)
\int_{M_{j}}|x| dx \\ &\leq\sum\limits_{j=1}\limits^{\infty}\phi(| \beta_{i_{j}}|)(c_{j}+r_{j})\lambda(M_{j})\\
&=\sum\limits_{j=1}\limits^{\infty}\phi(|\beta_{i_{j}}|)
\frac{1}{2^{i_{j}}\phi(|\beta_{i_{j}}|)}=1<\infty.
\end{aligned}
\end{equation*}
Thus, $h\in L_{\dagger}^{\phi}(\mu_{n})$.
Since $\phi$ satisfies the $\Delta_{2}$ condition, it follows that $h\in L^{\phi}(\mu_{n})$.
Observing that $\int_{M_{j}}xdx=\lambda(M_{j})(c_{j},0,\ldots,0)^{t}$ for $j\in\mathbb{N}$, we have
\begin{align*}
m(\xi\circ h)=\int_{\mathbb{R}^n}(\xi\circ h)(x)xdx&=\left(\sum\limits_{j=1}\limits^{\infty}\xi(\beta_{i_{j}})c_{j}\lambda(M_{j}),0,\ldots,0\right)^{t}.
\end{align*}
Additionally, from $\xi(\beta_{i_{j}})>2^{i_{j}}\phi(|\beta_{i_{j}}|)$ for $j\in\mathbb{N}$, we have
\begin{align*}
\sum\limits_{j=1}\limits^{\infty}\xi(\beta_{i_{j}})c_{j}\lambda(M_{j})
&>\sum\limits_{j=1}\limits^{\infty}2^{i_{j}}\phi(| \beta_{i_{j}}|)c_{j}\lambda(M_{j})\\
&=\sum\limits_{j=1}\limits^{\infty}\frac{c_{j}}{c_{j}+r_{j}}=\lim\limits_{l\rightarrow\infty}\sum\limits_{j=1}\limits^{l}\frac{c_{j}}{c_{j}+r_{j}}=\infty.
\end{align*}
Thus,
\begin{equation*}
\left|m(\xi\circ h)\right|=\infty.
\end{equation*}
If $\xi(\beta_{i_{j}})<-2^{i_{j}}\phi(|\beta_{i_{j}})|$ for $j\in\mathbb{N}$, a similar argument as in the first case yields
\begin{equation*}
\left|m(\xi\circ h)\right|=\infty,
\end{equation*}
which contradicts the fact that
\begin{equation*}
\left|m(\xi\circ h)\right|<\infty.
\end{equation*}
\qed

\begin{lem}\label{lem17}
Suppose that $\phi$ is a Young function satisfying the $\Delta_{2}$ condition. If an operator $\Psi: L^{\phi}(\mu_{n})\rightarrow\mathbb{R}^n$ is a continuous and $\mathrm{SL}(n)$ covariant valuation, then there is $\xi\in C_{\phi}(\mathbb{R}),$  and
\begin{equation*}
\Psi(h)=m(\xi\circ h)
\end{equation*}
for any $h\in L^{\phi}_{a}(\mu_{n}).$
\end{lem}
\emph{Proof.} We first prove that this result holds for  $s\in L^{\phi}_{a}(\mu_{n})$. Assume that $s\in L^{\phi}_{a}(\mu_{n})$, then there exists pairwise disjoint sets $M_{i}\in\mathscr{B}$ with $M_{i}\subset A(r,R]$ for $i=1,\ldots,l,$ and $\alpha_{1},\ldots,\alpha_{l}\in\mathbb{R}$ such that $s=\sum\limits_{i=1}\limits^{l}\alpha_{i}\chi_{M_{i}}$.
And Lemma \ref{lem12} and Lemma \ref{lem13} imply that
 \begin{align*}
 \Psi(s)&=\Psi\left(\sum\limits_{i=1}\limits^{l}\alpha_{i}\chi_{M_{i}}\right)=\sum\limits_{i=1}\limits^{l}\Psi(\alpha_{i}\chi_{M_{i}})\\
 &=\sum\limits_{i=1}\limits^{l}\xi(\alpha_{i})\int_{M_{i}}xdx
=\int_{\mathbb{R}^n}\left(\sum\limits_{i=1}\limits^{l}\xi(\alpha_{i})\chi_{M_{i}}(x)\right)xdx\\
 &=\int_{\mathbb{R}^n}(\xi\circ s)(x)x dx=m(\xi\circ s).
\end{align*}
From $\xi(0)=0$ and $\xi\circ s=\sum\limits_{i=1}\limits^{l}\xi(\alpha_{i})\chi_{M_{i}}$, we obtain the last equality.
\par 
Since simple functions are dense in $L^{\phi}_{a}(\mu_{n})$, for any $h\in L^{\phi}_{a}(\mu_{n}),$ there exists a sequence $\{s_{i}\}$ such that $\|s_{i}-h\|_{\phi} \rightarrow 0$. By Lemma \ref{lem6}, there is a subsequence $\{s_{i_{j}}\}$ of $s_{i}$ satisfying that $s_{i_{j}}(x) \rightarrow h(x)$ a.e. Furthermore, by the continuity of $\xi$, we obtain
\begin{equation*}
(\xi\circ s_{i_{j}})(x)x\rightarrow(\xi\circ h)(x)x, \,\,a. e.
\end{equation*}
Because $s_{i_{j}}\in L^{\phi}_{a}(\mu_{n})$, it follows that $|m(\xi\circ  s_{i_{j}})|<\infty$. By applying Lemma \ref{lem15}, we have
\begin{equation*}
|(\xi\circ s_{i_{j}})(x)x|\leq|(\xi\circ s_{i_{j}})(x)||x|\leq\lambda\phi(\vert s_{i_{j}}(x)\vert)|x|.
\end{equation*}
As $\lambda\phi(|s_{i_{j}}(x)|)|x|$ is integrable, by the dominated convergence theorem, we obtain
\begin{equation*}
\lim\limits_{j\rightarrow\infty}\int_{\mathbb{R}^n}(\xi\circ s_{i_{j}})(x)xdx=\int_{\mathbb{R}^n}(\xi\circ h)(x)xdx.
\end{equation*}
It is evident that $\|s_{i_{j}}-h\|_{\phi} \rightarrow 0$. By the continuity of $\Psi$, it follows that $\Psi(s_{i_{j}})\rightarrow\Psi(h)$ as $j\rightarrow\infty$.
This yields
\begin{equation*}
\lim\limits_{j\rightarrow\infty}\int_{\mathbb{R}^n}(\xi\circ s_{i_{j}})(x)xdx=\int_{\mathbb{R}^n}(\xi\circ h)(x)xdx=\Psi(h)
\end{equation*}
for any $h\in L^{\phi}_{a}(\mu_{n})$.
\begin{lem}\label{lem17}
Suppose that $\phi$ is a Young function satisfying the $\Delta_{2}$ condition. If an operator $\Psi: L^{\phi}(\mu_{n})\rightarrow\mathbb{R}^n$ is a continuous and $\mathrm{SL}(n)$ covariant valuation, then there is $\xi\in C_{\phi}(\mathbb{R}),$
and
\begin{equation}\label{001}
\Psi(h)=m(\xi\circ h)
\end{equation}
for every $h\in L^{\phi}(\mu_{n}).$
\end{lem}
\emph{Proof.}
 By Lemma \ref{lem4} and the valuation property of $\Psi$, we obtain
\begin{equation*}
\Psi(h)=\Psi(h)+\Psi(0)=\Psi(h\vee0)+\Psi(h\wedge0).
\end{equation*}
Therefore, it is enough to demonstrate that (\ref{001}) holds for all \( f \in L^{\phi}(\mu_{n}) \) for which \( h \geq 0 \) and for those for which \( h \leq 0 \).
Let $h\in L^{\phi}(\mu_{n})$ with $h\geq 0$. For $k\in\mathbb{N}$, we define $h_{k}:\mathbb{R}^n\rightarrow \mathbb{R}$ as
\begin{equation*}
h_{k}=\sum\limits_{j=0}\limits^{k}\chi_{M_{j}}h,
\end{equation*}
where $M_{j}=A[1/2^{j+1},1/2^{j})\cup A[j+1,j+2)$ for $j\geq0$.
This leads that $0\leq h_{1}\leq \cdots h_{k}\leq h$ and $h_{k}\rightarrow h$, a.e. Due to the continuity of $\xi$ , it follows that
\begin{equation*}
(\xi\circ h_{k})(x)x\rightarrow(\xi\circ h)(x)x, \,\,\mathrm{a.e.}
\end{equation*}
Because $h_{k}\in L^{\phi}_{a}(\mu_{n})$, it implies $|\Psi(h_{k})|<\infty$. By applying Lemma \ref{lem15}, we obtain
\begin{equation*}
|(\xi\circ h_{k})(x)x|\leq|(\xi\circ h_{k})(x)|| x|\leq\lambda\phi(h_{k}(x))|x|\leq\lambda\phi(h(x))|x|, \,\,\mathrm{a.e.}
\end{equation*}
As $\lambda\phi(h(x))|x|$ is integrable, by the dominated convergence theorem, we have
\begin{equation*}
\lim\limits_{k\rightarrow\infty}\int_{\mathbb{R}^n}(\xi\circ h_{k})(x)xdx=\int_{\mathbb{R}^n}(\xi\circ h)(x)xdx.
\end{equation*}
Moreover, since $|h_{k}-h|\leq h$ and $\phi$ is increasing function, we have $\phi(|h_{k}-h|)\leq\phi( h)$. As $\phi( h)$ is integrable, by the dominated convergence theorem, we deduce that
$\rho_{\phi}(|h_{k}-h|)\rightarrow 0$, which leads to  $\|h_{k}-h\|_{\phi}\rightarrow 0$. Due to the continuity of $\Psi$, we have $\Psi(h_{k})\rightarrow\Psi(h)$ as $k\rightarrow\infty$. Therefore, it concludes that
\begin{equation*}
\Psi(h)=\int_{\mathbb{R}^n}(\xi\circ (h)(x)xdx.
\end{equation*}
Let $h\in L^{\phi}(\mu_{n})$ with $h\leq 0$, we similarly define
\begin{equation*}
h_{k}=\sum\limits_{j=0}\limits^{k}\chi_{M_{j}}h,
\end{equation*}
where $M_{j}=A[1/2^{j+1},1/2^{j})\cup A[j+1,j+2)$ for $j\geq0$.
This implies that $h \leq h_{k}\leq \cdots h_{1}\leq 0$ and $h_{k}\rightarrow h$, a.e. Following  the same argument, we have
\begin{equation*}
\Psi(h)=\int_{\mathbb{R}^n}(\xi\circ h(x))x dx.
\end{equation*}

\qed

\par \textbf{Acknowledgement}\quad The first author expresses her sincere gratitude to Monika Ludwig for her guidance, support, and encouragement.





\begin{thebibliography}{99}

\bibitem{Alesker1} S. Alesker, Continuous rotation invariant valuations on convex sets, Ann. of Math. (2) 149 (1999), 977-1005.

\bibitem{Alesker2} S. Alesker, Description of continuous isometry covariant valuations on convex sets, Geom. Dedicata 74 (1999), 241-248.

\bibitem{Alesker3} S. Alesker, On P. McMullen's conjecture on translation invariant valuations, Adv. Math. 155 (2000), 239-263.

\bibitem{Alesker4} S. Alesker, Description of translation invariant valuations on convex sets with solution of P. McMullen's conjecture, Geom. Funct. Anal. 11 (2001), 244-272.





\bibitem{Colesanti1} A. Colesanti, D. Pagnini, P. Tradacete, and I. Villanueva, A class of invariant valuations on $\mathrm{Lip}(S^{n-1})$, Adv. Math. 366 (2020): 107069.

\bibitem{Colesanti2} A. Colesanti, D. Pagnini, P. Tradacete, and I. Villanueva, Continuous valuations on the space of Lipschitz functions on the sphere, J. Funct. Anal. 280 (2021): 108873.



\bibitem{Klain1} D. A. Klain, Star valuations and dual mixed volumes, Adv. Math. 121 (1996), 80-101.

\bibitem{Klain2} D. A. Klain, Invariant valuations on star-shaped sets, Adv. Math. 125 (1997), 95-113.

\bibitem{Klain3} D. A. Klain, Even valuations on convex bodies, Trans. Amer. Math. Soc. 352 (2000), 71-93.

\bibitem{Klain-Rota} D. A. Klain, G. C. Rota, Introduction to geometric probability,  Lezioni Lincee [Lincei Lectures], Cambridge University Press, Cambridge, 1997.


\bibitem{Kone} H. Kone, valuations on Orlicz spaces and $L^{\phi}$-star sets, Adv. in Appl. Math. 52 (2014) 82-98.



\bibitem{Ludwig1} M. Ludwig, Moment vectors of polytopes, IV International Conference in "Stochastic Geometry, Convex Bodies, Empirical Measures and Applications to Engineering Science", Vol. II (Tropea, 2001), Rend. Circ. Mat. Palermo (2) Suppl. 70 (2002), 123-138.

\bibitem{Ludwig2} M. Ludwig, Projection bodies and valuations, Adv. Math. 172 (2002), 158-168.

\bibitem{Ludwig3} M. Ludwig, Valuations of polytopes containing the origin in their interiors, Adv. Math. 170 (2002), 239-256.

\bibitem{Ludwig4} M. Ludwig, Ellipsoids and matrix-valued valuations, Duke Math. J. 119 (2003), 159-188.

\bibitem{Ludwig9} M. Ludwig, Minkowski valuations, Trans. Amer. Math. Soc. 357 (10) (2005), 4191-4213.

\bibitem{Ludwig5} M. Ludwig, Valuations on Sobolev spaces, Amer. J. Math. 134 (3) (2012), 827-842.

\bibitem{Ludwig6} M. Ludwig, Covariance matrices and valuations, Adv. in Appl. Math. 51 (3) (2013), 359-366.

\bibitem{Ludwig7} M. Ludwig, Fisher information and matrix-valued valuations, Adv. Math. 226 (3) (2011), 2700-2711.

Fond. CIME/CIME Found. Subser.
Springer, Cham, 2023, 19–78.


\bibitem{Ludwig8} M. Ludwig, M. Reitzner, A characterization of affine surface area, Adv. Math. 147 (1999), 138-172.

\bibitem{L-Y-Z1} E. Lutwak, D. Yang, and G. Zhang, The Cramer-Rao inequality for star bodies, Duke Math. J. 112 (2002), 59-81.



\bibitem{McMullen1} P. McMullen, Valuations and dissections, Handbook of convex geometry, Vol. A, B, North-Holland Publishing Co., Amsterdam, 1993, 933-988.

\bibitem{McMullen2} P. McMullen, R. Schneider, Valuations on convex bodies, Convexity and its applications, Birkh\"{a}user Verlag, Basel, 1983, 170-247.

\bibitem{Ma} D. Ma, Moment matrices and  $\textrm{SL}(n)$ equivariant valuations on polytopes, Int. Math. Res. Not. 14 (2021), 10436-10489.


\bibitem{Parapatits1} L. Parapatits, $\textrm{SL}(n)$-contravariant $L_p$-Minkowski valuations, Trans. Amer. Math. Soc. 366 (2014), 1195-1211.

\bibitem{Parapatits2} L. Parapatits, $\textrm{SL}(n)$-covariant $L_p$-Minkowski valuations, J. Lond. Math. Soc. 89 (2) (2014), 397-414.



\bibitem{Tsang1} A. Tsang, Valuations on $L^{p}$ spaces, Int. Math. Res. Not. 20 (2010), 3993-4023.

\bibitem{Tsang2} A. Tsang, Minkowski valuations on $L^{p}$-spaces, Trans. Amer. Math. Soc. 364 (2012), 6159-6186.

\bibitem{Wannerer} T. Wannerer, $\textrm{GL}(n)$ equivariant Minkowski valuations, Indiana Univ. Math. J. 60 (2011), 1655-1672.

\bibitem{Wang} T. Wang, Semi-valuations on $\mathrm{BV}(\mathbb{R}^n)$, Indiana Univ. Math. J. 63 (5) (2014), 1447-1465.


\bibitem{Wang1} W. Wang, $\textrm{SL}(n)$ equivariant matrix-valued valuations on polytopes, Int. Math. Res. Not. 13 (2022), 10302-10346.

\bibitem{Wang3} W. Wang, R. G. He, and L. J. Liu, $\textrm{SL}(n)$ covariant vector-valued valuations on $L^{p}$-spaces, Ann. math. du. Qu$\acute{e}$. 45 (2021), 465-486.




\bibitem{Zeng-Ma} C. N. Zeng, D. Ma, $\textrm{SL}(n)$ covariant vector valuations on polytopes, Trans. Amer. Math. Soc. 370 (2018), 8999-9023.

\bibitem{Zeng-Zhou} C. N. Zeng, Y. Q. Zhou, $\textrm{SL}(n)$ contravariant matrix-valued valuations on polytopes, Int. Math. Res. Not. 15(2024), 11343–11373.



\end{thebibliography}
\end{document}